\documentclass[graybox]{svmult}

\usepackage{type1cm}        
\usepackage{makeidx}         
\usepackage{graphicx}        
\usepackage{multicol}        
\usepackage[bottom]{footmisc}

\usepackage{newtxtext}       %
\usepackage{newtxmath}       

\usepackage[normalem]{ulem} 

\makeindex             
\newtheorem{prn}{Proposition}
\newtheorem{thm}{Theorem} 
 
\newcommand{\supp}{\mathop{\mathrm{supp}}}

\newcommand{\Trunc}{\mathop{\mathrm{Trunc}}}
\newcommand{\trunc}{\mathop{\mathrm{trunc}}}
\def\bfl{\mbox{\boldmath$\lambda$}}
\newcommand{\R}{\mathbb{R}}
\def\bfz{\mbox{\boldmath$\zeta$}}

\begin{document}

\title*{THB-spline approximations for turbine blade design with local B-spline approximations}
%
\author{\centerline{Cesare Bracco, Carlotta Giannelli, David Gro{\ss}mann, Sofia Imperatore,}  \centerline{Dominik Mokri\v{s}, Alessandra Sestini}}
\authorrunning{C. Bracco et. al.} 
\institute{
Cesare Bracco, Carlotta Giannelli, Sofia Imperatore, Alessandra Sestini\at Dipartimento di Matematica e Informatica ``U. Dini'', Universit\`a degli Studi di Firenze, Viale Morgagni  67/A, 50134 Firenze, Italy \email{(cesare.bracco, carlotta.giannelli, alessandra.sestini)@unifi.it, sofia.imperatore@stud.unifi.it}
\and David Gro{\ss}mann, Dominik Mokri\v{s}\at MTU Aero Engines AG, Munich, Germany \email{(David.GROSSMANN, Dominik.Mokris)@mtu.de}
}
%
%
\maketitle
\abstract*{
} 
\abstract{
We consider two-stage scattered data fitting with truncated hierarchical B-splines (THB-splines) for the adaptive reconstruction of industrial models. The first stage of the scheme is devoted to the computation of local  least squares variational spline approximations, exploiting a simple fairness functional  to handle data distributions with a locally varying density of points. Hierarchical spline quasi-interpolation based on THB-splines is considered in the second stage of the method to construct the adaptive spline surface approximating the whole scattered data set and a suitable strategy to guide the adaptive refinement is introduced. 
A selection of examples on geometric models representing components of aircraft turbine blades highlights the performances of the scheme. The tests include a scattered data set with voids and the adaptive reconstruction of a cylinder-like surface.}

\section{Introduction}\label{sec:1}

Scattered data fitting is nowadays of fundamental importance in a variety of fields, ranging from geographic applications to medical imaging and geometric modeling. 
The topic can be addressed in different approximation spaces, either by using  spline spaces or radial basis functions. In particular, thanks to their low computational cost and also to their better control of conditioning, schemes relying on local approximations  have received a lot of attention, being formulated either as partition of unity interpolation or as two-stage approximation methods, see for example \cite{bracco2017,CDP18,DMS06,DS08} and references  therein. 

In this work we are interested in  two-stage surface reconstruction of industrial models starting from scattered data obtained by optical scanning acquisitions. In order to reduce the noise, the available data are always preprocessed because the interest is in highly accurate reconstructions.  Anyway, since the available data can not be considered  exact, interpolation is not required and a less costly approximation scheme can be used to compute any local fitting. Furthermore, considering that splines represented in B-spline form are the standard choice for industrial computer-aided design applications (see e.g. \cite{prautzsch2002}), we are interested in determining the final approximation in a space spanned by suitable extensions of tensor-product B-splines.

As well known, an approximating spline can be obtained by using different approximation approaches and also operating in different kinds of spline spaces, with the related schemes usually  divided into two main classes. The first is the class of \emph{global} schemes which use simultaneously all the given information and thus require the solution of a linear system whose size is  (substantially) equal to the cardinality of the considered spline space.  The second class collects all kinds of \emph{local} schemes which avoid the solution of a global linear system. Intermediate alternatives are also possible, see e.g., \cite{LaiSchum09}, where spline spaces on triangulations are considered together with a  domain decomposition technique. Many local methods can be collocated in the field of  Quasi-Interpolation (QI) which is a fundamental methodology in the context of spline approximation, see e.g. \cite{Lee01,LycheSchum75} for an introduction. In some cases, e.g. when dealing with scattered data, QI can require the solution of small linear systems whose size does not depend on the cardinality of the entire data set or on the size of the global spline space. These local systems descend from the application of \emph{local} approximation schemes, each one considering a small number of data whose associated parameter values belong to a  modest portion of the parametric domain intersecting the support of a  compactly supported basis function. The computational advantage of the QI approach is evident, since these local linear systems are  of small size and independent of each other. However, the development of an effective scattered data fitting approach  based on local methods is never trivial, since  the quality of the final spline approximation depends in this case not only on the considered global spline space but also on the approximation power of the local  scheme and on the choice of the  corresponding local data set.
Here we are interested in considering a method of this kind but having also the adaptivity feature, in particular relying on a QI two stage methodology and working in adaptive spline spaces, since they allow local refinement and generalize the standard tensor-product model. Since  scattered data can be characterized by a highly varying distribution, by also including voids, a flexible and reliable approach which automatically (re-) constructs the geometric model may strongly improve the efficiency of the overall scheme by suitably adapting the solution to the shape and configuration of the given point clouds. 
 
 In recent years the need for local refinement both in computer aided geometric design and isogeometric analysis has motivated a lot of research on adaptive spline spaces, leading to the introduction of several  adaptive spline constructions, such as spline spaces over T-meshes \cite{deng2006}, T-splines \cite{sederberg03}, LR B-splines \cite{dokken} and hierarchical B-splines \cite{kraft}. Spline spaces over T-meshes is the most straightforward approach: given a T-mesh, that is a rectangular mesh allowing T-junctions (vertices with valence $3$), a space of functions which are polynomial over each cell can be easily defined. While this is very natural, it has been shown that the dimension of the corresponding space is stable (namely, it depends only on the topology of the mesh and not on the geometry) only under certain conditions, see \cite{mourrain} and \cite{bracco2016b}.  Consequently, an efficient construction of a global basis is, in general, still an open question. T-splines \cite{sederberg03} and LR B-splines \cite{dokken} are both formally B-splines, but defined on local knot vectors depending on the topology of the mesh, and therefore allowing local refinement. They are very flexible, but their linear independence is guaranteed only under certain conditions on the mesh (different for the two types of spline), which any refinement algorithm must then preserve, see e.g. \cite{daveiga} and \cite{patrizi}. The construction of hierarchical B-splines \cite{kraft} is based on a multi-level construction, where the refinement is obtained by replacing coarser elements with finer ones in the mesh hierarchy. As a consequence, the corresponding basis is composed of B-splines constructed on meshes of different levels. The multilevel nature of these splines allows the design of efficient local refinement schemes and their suitable integration into existing computer aided design software \cite{kiss2014}. Moreover, the truncated version of hierarchical B-splines (THB-splines) \cite{giannelli2012} can be used for an easy extension to hierarchical spline spaces of any QI scheme formulated in a standard spline space \cite{speleers2017,speleers2016}. These features make THB-splines a natural choice to solve our scattered data fitting problem and are then the solution here considered.

The first proposal based on adaptive THB-spline fitting of scattered data for the reconstruction of industrial models was introduced in \cite{kiss2014}, where a global adaptive smoothed least-squares scheme was developed. Successively, in order to increase its locality and reduce the computational cost, the same problem has been addressed  in different papers by combining a two-stage approach with hierarchical spline approximations. The first contribution where this kind of schemes was used by some of the authors was presented in \cite{bracco2016}  to deal with gridded data of Hermite type.  In \cite{bracco2017} these kinds of approximants were extended for the first time to scattered data, by using in the first stage of the scheme local polynomial least squares approximations of variable degree. A preliminary application of this scheme to industrial data reconstruction was given in \cite{bracco2018}, where its theoretical analysis was also presented. 

In this paper a variant of the approach considered in \cite{bracco2017,bracco2018} is presented,  to further decrease the number of degrees of freedom necessary to reach a certain accuracy and also to reduce the artifacts in the resulting surface. Since the local polynomials used in \cite{bracco2017,bracco2018} need to be converted into B-spline form for being usable by the quasi-interpolation approach considered in the second stage of scheme, we now work directly in local spline spaces.  In this way, the algorithm for the first stage of the scheme is  simplified. In order to avoid rank-deficiency problems by simultaneously controlling the smoothness of local approximants, a smoothing term is added to the least squares objective function.   To improve the stability of the proposed method for general data configuration, we also inserted an automatic control on the choice of the local data sets. 
 
The structure of the paper is as follows. Section~\ref{sec:1bis} presents the model problem, while the first stage of the new scheme, devoted to the computation of local smoothed least squares B-spline approximations, is described in Section~\ref{sec:2}. Section~\ref{sec:3} introduces the construction of the adaptive THB-spline surface approximating the scattered data set. A selection of examples on geometric models representing components of aircraft turbine blades is presented in Section~\ref{sec:4}, and compared with the results obtained in \cite{bracco2018}. The numerical experiments include a new scattered data set with voids and the adaptive reconstruction of  a surface closed in one parametric direction.

\section{The problem} \label{sec:1bis}
The industrial models here considered are components of aircraft turbine blades which can be  suitably represented in parametric form by using just one map,  with the possibility of being periodic in one parametric direction. The problem  can be mathematically described as follows. Let
\[
F := \left\{\mathbf{f}_i \in\mathbb{R}^3, i=1,\ldots,n\right\}
\]
be a scattered set of distinct points  in the 3D space which can be reasonably associated by a one-to-one map to a set $X := \left\{\mathbf{x}_i := (x_i,y_i) \in \Omega \subset\mathbb{R}^2, i=1,\ldots,n\right\}$ of
distinct parameter values belonging to a closed planar parametric domain $\Omega.$   Since the choice of a suitable parametrization method for the definition of the set $ X$, which can naturally influence the quality of the final approximation,  is not our focus, in this paper we relate to classical parametrization methods based on a preliminary triangulation of the scattered data set $F,$ see e. g., \cite{floater1997,floater2005}. Consequently, both $F$ and $X$ are considered input data for the approximation problem. 

Focusing on two-stage spline approximation schemes, we can introduce the general idea referring for simplicity to their formulation in a standard space $V$  of tensor product splines of bi-degree  $\mathbf{d} = (d_1,d_2),$ where it is assumed $d_i \ge 1, i=1,2,$ in order to deal with at least continuous functions. With this kind of methods, a quasi-interpolation operator $Q$ is defined so that $Q(F,X) = \mathbf{s},$ with $\mathbf{s}$ denoting a vector function, possibly periodic in one parametric direction, with components in the spline space $V.$ Using  a suitable spline basis ${\cal B} := \left\{B_J \right\}_{J \in \Gamma}$ of $V,$ such vector spline $\mathbf{s}$ can be expressed as follows,
\begin{equation}\label{eq:qi}
\mathbf{s}(\mathbf{x}) := \sum_{J \in \Gamma} \bfl_J(F_J,X_J) B_J(\mathbf{x})\,, \qquad \mathbf{x} \in  {\Omega}, \end{equation}
where each coefficient vector $\bfl_J(F_J,X_J), \, J \in \Gamma,$ is computed in the first stage of the scheme by using a certain local subset $F_J \subset F$ of data and the corresponding set of parameter values $X_J \subset X$, so that $\mathbf{s}(\mathbf{x}_i) \approx \mathbf{f_i}, \, i=1,\ldots,n.$ 

When dealing with discrete data, measuring the accuracy of the spline approximation with the maximum  of the errors $||\mathbf{s}(\mathbf{x}_i) - \mathbf{f_i}||_2$ at each parameter site can appear reasonable at the first sight. However, the quality of the approximant is also strictly related to the lack of unwanted artifacts,  a feature of fundamental importance for industrial applications of any approximation scheme. In this context, it is then a common practice to require the error under a prescribed tolerance only at a certain percentage of sites in $X$.



\section{First stage: local B-spline approximations}
\label{sec:2}
For computing each vector coefficient $\bfl_J,\, J\in \Gamma,$ necessary in  \eqref{eq:qi} to define  the approximation $\mathbf{s},$ we consider a local data subset $F_J\subset F$,
 \[
F_J := \left\{\mathbf{f}_i: i\in I_J\right\} \quad 
\text{with} \quad I_J := \left\{ i:  \mathbf{x}_i \in X \cap \Omega_J\right\}, 
\]
associated to the set 
\[
X_J := \{\mathbf{x}_i: i \in  I_J\} 
\]
of parameter values in a local subdomain $\Omega_J$ of $\Omega$ which has non empty intersection with the support of the basis function $B_J$, namely
$\Omega_J \cap \supp(B_J) \ne \emptyset$. By denoting with
${\cal B}_J := \left\{B_I : I \in \Lambda_J \subset \Gamma \right\}$ the set of B-splines  in ${\cal B}$  not vanishing in $\Omega_J$ (which necessarily includes $B_J$), the value of $\bfl_J $ is defined as the vector coefficient associated with $B_J$ in a local spline approximation  $\mathbf{s}_J \in {\cal S}_J,$ with $\mathbf{s}_J: \Omega_J \rightarrow \R^3$ and 
$${\cal S}_J := \mathrm{span}\{B_I, I \in \Lambda_J\}.$$
Concerning this local spline space,  the following proposition is proved, since it shows that ${\cal S}_J$ has reasonable approximation power and in particular that it includes the restriction to $\Omega_J$ of any linear polynomial,
 \begin{prn}
{\it The following space inclusion holds true}
$$\Pi^2_{\bf d} \vert_{\Omega_J} \subseteq {\cal S}_J:= \mathrm{span}\{B_I : I \in \Lambda_J \}\,,$$
{\it where} $\Pi^2_{\bf d} \vert_{\Omega_J}$ {\it denotes the restriction to} $\Omega_J$ {\it of the tensor product space of bivariate polynomials of bi-degree} ${\bf d}.$
 \begin{proof}
Let $c_K$ be a cell of ${\cal G}$ such that $c_K \cap \Omega_J \ne \emptyset.$ Then the definition of $\Lambda_J$ implies that  $c_K \subseteq \mathrm{supp}(B_I) \, \Rightarrow \, I \in \Lambda_J.$ Thus we can say that $\mathrm{span}\{ B_I : c_K \subseteq \mathrm{supp}(B_I)\} \subseteq {\cal S}_J.$ 
Since $\mathrm{span}\{ B_I : c_K \subseteq \mathrm{supp}(B_I)\} = \Pi^2_{\bf d} \vert{c_K},$ the proof is completed considering that $c_K$ is any cell of ${\cal G}$ with non vanishing intersection with $\Omega_J$.
\end{proof}
\end{prn}

\noindent
 The variational fitting method adopted to determine $\mathbf{s}_J$ in ${\cal S}_J$ consists in minimizing the following objective function
\begin{equation} \label{obj}
\sum_{i \in I_J} 
\Vert \mathbf{s}_J(\mathbf{x}_i) - \mathbf{f}_i \Vert_2^2 + \mu \, E(\mathbf{s}_J)\,,
\end{equation}
where $\mu > 0 $ is a smoothing coefficient and $E(s_J)$ the thin-plate energy,
\[
E(\mathbf{s}_J) := \int_{\Omega_J} 
\left|\left|\frac{\partial^2 \mathbf{s}_J}{\partial x^2}\right|\right|_2^2
+ 2 \left|\left|\frac{\partial^2  \mathbf{s}_J}{\partial x \partial y}\right|\right|_2^2
+\left|\left|\frac{\partial^2  \mathbf{s}_J}{\partial y^2}\right|\right|_2^2 \,\text{d}x\text{d}y\,.
\]

\noindent
As recalled in the Appendix, the assumption of a positive $\mu$ ensures that this local approximation problem admits always a unique solution,  provided that the sites belonging to $X_J$ are not collinear.   

 Since the scheme is locally applied, an automatic (data-dependent) selection of the parameter $\mu$ could be considered. For example, the choice may take into account the cardinality $|X_J|$ of the local sample or the area of $\Omega_J$, which influence the value of the first and of the second addend in (\ref{obj}), respectively.  In view of this influence however, we may observe that a constant value of $\mu$ implies that the balancing between the fitting and the smoothing term in the objective function usually increases when $|X_J|$ or the area of $\Omega_J$ increases, being this true in the second case because second derivatives are involved in the smoothing term. Both these behaviors seem reasonable and are confirmed by the quality of the results obtained in our experiments, where a constant value for $\mu$ is suitably chosen.

Differently from \cite{bracco2017,bracco2018}, in order to better avoid overfitting, a lower bound $n_{\mathrm{min}}$ for the cardinality of $X_J$ is now required, being $n_{\mathrm{min}} \ge 3$ the only additional input parameter required by the algorithm, besides $\mu.$  To fulfill this condition, $\Omega_J$ is initialized as $\supp(B_J)$ and enlarged until $|X_J| \ge n_{\mathrm{min}}$. Note that the refinement strategy presented in Section \ref{sec:3} automatically guarantees that the inequality $|X_J| \ge n_{\mathrm{min}}$ becomes fulfilled without requiring  an excessive enlargement of the set $\Omega_J$, which would compromise the locality of the approximation. For this reason, it is not necessary to set a maximum value for controlling the enlargement of the local data set.  Considering that the smaller $n_{\mathrm{min}}$ is, the higher is the obtainable level of detail but also the probability of overfitting, in our experiments (which always adopt $d_1 = d_2 = d = 2, 3$), a good low range for its selection has always been $ d^2 \le n_{\mathrm{min}} \le (d+1)^2.$ 

\noindent
Summarizing, the computation of each $\bfl_J,\, J\in \Gamma,$ is done according to the following algorithm.

\medskip
\noindent {\bf Algorithm 1: local smoothing spline approximant}

\medskip
\noindent {\bf Inputs}
\begin{itemize}
\item $F\subset \mathbb{R}^3$: scattered data set;
\item $X\subset \Omega \subset \mathbb{R}^2$: set of parameter values corresponding to the data in $F$;
\item ${\cal G}$: uniform tensor-product mesh in $\Omega$ (possibly with auxiliary cells); 
\item $V$: tensor-product spline space of bi-degree ${\bf d}$ associated with ${\cal G}$ with B-spline basis ${\cal B}$; 
\item $J \in \Gamma$: index of the considered basis function $B_J\in{\cal B}$;
\item $\mu$: smoothing spline parameter ($\mu > 0$);
\item $n_{\mathrm{min}} $: minimum required number of local data ($3 \le n_{\mathrm{min}} \ll \vert F \vert$);
\end{itemize}
 
\begin{enumerate}
\item Initialization
\begin{enumerate}
\item initialize $\Omega_J=\supp(B_J)$;
\item initialize  $I_J=\left\{ i:  \mathbf{x}_i \in X \cap \Omega_J\right\}, \, F_J=\left\{\mathbf{f}_i: i \in I_J\right\}$ and $X_J=\left\{\mathbf{x}_i: i \in I_J\right\}$; 
\end{enumerate}
\item while  $\vert F_J \vert < n_{\mathrm{min}}$   
\begin{enumerate}
       \item  enlarge $\Omega_J$ with the first ring of cells in ${\cal G}$ surrounding $\Omega_J$;
       \item update $I_J=\left\{ i:  \mathbf{x}_i \in X \cap \Omega_J\right\}, \, F_J=\left\{\mathbf{f}_i: i \in I_J\right\}$ and $X_J=\left\{\mathbf{x}_i: i \in I_J\right\}$;
       \end{enumerate}
\item if the sites in $X_J$ are not collinear, then:
\begin{enumerate}
\item compute the local approximation $\mathbf{s}_J = \sum_{I \in \Lambda_J} \mathbf{c}_I^{(J)} B_I$  minimizing the objective function in \eqref{obj} for the data $F_J$ and $X_J$;
\item set $\bfl_J = \mathbf{c}_J^{(J)};$
\end{enumerate}
else set $\bfl_J = \frac{1}{ \vert F_J \vert} \ \displaystyle{\sum_{I \in I_J}} \mathbf{f}_I.$ 
\end{enumerate} 

\noindent{\bf Output} 
\begin{itemize}
\item  vector coefficient $\bfl_J$.
\end{itemize}

Note that  exploiting a regularized least square approximation and, as a consequence, being able to directly employ the local spline space has significantly simplified  the   algorithm originally proposed for the first stage in \cite{bracco2017,bracco2018}, where a variable-degree local polynomial approximation was considered. In particular,  the new scheme does not require the selection of a suitable  degree for the computation of any coefficient $\bfl_J$ and eliminates the conversion of the computed approximant from the polynomial to the B-spline basis.

In the following section, after introducing the THB-spline basis, the operator $Q$ is easily extended to hierarchical spline spaces, by also introducing the automatic refinement algorithm here considered. Note that this extension rule ensures that the coefficient associated to a THB-spline basis function remains unchanged  on a refined hierarchical mesh if this function remains active on the updated hierarchical configuration.




\section{Second stage: THB-spline approximation}
\label{sec:3}
Let us consider a sequence $V^0\subset \ldots\subset V^{M-1}$ of  $M$ spaces of tensor-product splines of bi-degree $\mathbf{d}:=(d_1,d_2),$  $d_i \ge 1, i=1,2,$ defined on the closed domain $\Omega$, and each one associated with the tensor-product grid ${\cal G}^{\ell}$ and the basis ${\cal B}_{\mathbf d}^\ell := \left\{B_J^\ell\right\}_{J \in \Gamma_{\mathbf d}^\ell}$. Let $\Omega^0\supseteq \ldots \supseteq \Omega^{M}$ be a sequence of closed domains, with $\Omega^0 := \Omega$ and $\Omega^M := \emptyset$. Each $\Omega^\ell$, for $\ell=1,\ldots,M-1$ is the union of cells of the tensor-product grid ${\cal G}^{\ell-1}$. Let ${\cal G}_{\cal H}$ be the hierarchical mesh defined by
\begin{equation*}
{\cal G}_{\cal H} : = \{Q\in {\mathcal G}^\ell,\ 0\le \ell \le M-1\} \quad \mbox{with} \quad {\mathcal G}^\ell := \{Q\in G^\ell:~ Q\subset \Omega^\ell \backslash \Omega^{\ell+1} \},
\end{equation*}
where each ${\mathcal G}^\ell$ is called the set of {\it active cells} of level $\ell$.
The hierarchical B-spline (HB-spline) basis ${\mathcal H}({\cal G}_{\cal H})$ with respect to the mesh ${\cal G}_{\cal H}$ is defined as
\begin{equation*}
{\mathcal H}({\cal G}_{\cal H}):=\{B_J^\ell:\,J\in A^\ell, \ell=0,\ldots,M-1\},
\end{equation*}
where
\begin{equation*}
A^\ell := \left\{J\in \Gamma^\ell : \supp(B_J^\ell) \subseteq \Omega^\ell \wedge
\supp(B_J^\ell) \not\subseteq \Omega^{\ell+1}\right\},
\end{equation*}
is the set of \emph{active multi-indices} of level $\ell$, and $\supp(B_J^\ell)$ denotes the intersection of the support of $B_J^\ell$ with $\Omega^0$. The corresponding hierarchical space is defined as $S_{\cal H} := {\rm span}\ {\mathcal H}({\cal G}_{\cal H})$.

For any $ s \in V^\ell$, $\ell=0,\ldots,M-2$, let
\begin{equation*}
 s = \sum_{J\in \Gamma^{\ell+1}} \sigma_J^{\ell+1} B_J^{\ell+1}
\end{equation*}
be its representation in terms of B-splines of the refined space $V^{\ell+1}$.  We define the truncation of $s$ with respect to level $\ell+1$ and its (cumulative) truncation with respect to all finer levels as
\begin{equation*}
{\trunc}^{\ell+1}  (s) := \sum_{J\in \Gamma^{\ell+1}:\, 
\supp(B_J^{\ell+1})\not\subseteq\Omega^{\ell+1}} \sigma_J^{\ell+1} B_J^{\ell+1},
\end{equation*}
and
\begin{equation*}
{\Trunc}^{\ell+1}( s) := {\trunc}^{M-1}({\trunc}^{M-2}(\ldots({\trunc}^{\ell+1}(s))\ldots)), 
\end{equation*}
respectively. For convenience, we also define ${\Trunc}^{M}(s) := s$ for $s\in V^{M-1}$.
The THB-spline basis ${\cal T}({\cal G}_{{\cal H}})$ of the hierarchical space $S_{\cal H}$ was introduced in \cite{giannelli2012} and can be defined as 
\begin{equation*}
{\cal T}({\cal G}_{{\cal H}}) 
:=
\left\{
T_J^\ell:= {\Trunc}^{\ell+1} (B_J^\ell) : J\in A^\ell, \ell=0,\ldots,M-1
\right\}.
\end{equation*}
The B-spline $B_J^\ell$ is called the \emph{mother} B-spline of the truncated basis function $T_J^\ell$. 

 We recall that THB-splines are linearly independent, non-negative, preserve the coefficients of the underlying sequence of B-splines, and form a partition of unity \cite{giannelli2012,giannelli2014}. Besides that, following the general approach introduced in \cite{speleers2016}, using such basis we can easily construct the vector THB-spline approximation of the whole scattered data set in terms of the hierarchical quasi-interpolant  $\mathbf{s} = Q(F,X)$ as follows, 
\begin{equation}\label{hqi1}
\mathbf{s} (\mathbf{x}) := \sum_{\ell =0}^{M-1}  \sum_{I \in A^\ell} \bfl_I^\ell(F_I,X_I) T_I^\ell (\mathbf{x}),
\end{equation}
where each vector coefficient $\bfl_I^\ell$ is the one of the mother function $B_I^\ell$ and is obtained by computing the local regularized B-spline approximation  $\mathbf{s}_I^\ell$ on the data set $F_I^\ell$ associated to $B_I^\ell$ as described in Section~\ref{sec:2}. 
 
  In order to define an adaptive approximation scheme, the following algorithm is used to iteratively construct the hierarchical mesh ${\cal G_H},$  the corresponding spline space $S_{\cal H}$ and the final approximating spline $\mathbf{s}.$ As any automatic refinement strategy, the algorithm requires in input  some parameters which drive the refinement process. One of them is the error tolerance $\epsilon >0$ whose value has to be chosen not only considering the accuracy desired for the reconstruction but also taking into account the level of noise affecting the given points cloud (here assumed without significant outliers). The  percentage bound parameter $\eta$ specifies the number of data points for which the error is required to be within the given tolerance and a value strictly less than $100\%$ is suggested to reduce the influence of outliers with moderate size on the approximation. Another required parameter is the maximum number of levels $M_{\mathrm{max}}$ which has to be chosen considering the maximal   level of detail desired for the reconstruction. The choice of the other input parameters, $n_{\mathrm{loc}}, n_1$ and $n_2,$ is discussed after the algorithm. 
\bigskip

\noindent {\bf Algorithm 2: adaptive hierarchical spline fitting}
\medskip

\medskip
\noindent {\bf Inputs}
\begin{itemize}
\item $F \subset \mathbb{R}^3$: scattered data set with $n = |F|$;
\item $X \subset \Omega \subset \mathbb{R}^2$: set of parameter values corresponding to the data in $F$;
\item ${\cal G}^0$: initial uniform tensor-product mesh in $\Omega$ (possibly with auxiliary cells); 
\item $V^0 $: tensor-product spline space of bi-degree ${\bf d}$ associated with ${\cal G}^0$;
\item $\epsilon :$ error tolerance;
\item $\eta$: percentage bound of data points for which the error is required to be within the tolerance $\epsilon$ 
(default: $\eta = 95 \%$);
\item $M_{\mathrm{max}}:$  integer specifying the maximum number of levels;
\item $ n_{\mathrm{loc}}$: minimum number of local data required for refinement ($3 \le n_{\mathrm{loc}} \ll \vert F \vert$);
\item $n_1\,, n_2$: positive integers specifying the number of uniform horizontal or vertical splittings for the support of a tensor-product B-spline
(default: $n_1=n_2 = 1$). 
\end{itemize}

\begin{enumerate}
\item Initialization
\begin{enumerate}
\item set ${\cal G_H}={\cal G}^0$and  $S_{\cal H}=V^0;$ 
\item set the current number of level $M=1;$ 
\item use Algorithm 1 to compute the coefficients of the  hierarchical QI vector spline $\mathbf{s} \in S_{\cal H}$ introduced \eqref{hqi1}. 
\item evaluate the errors at the data sites $\mathbf{x}_i, \in X, i=1,\ldots,n\,,$
\begin{equation*}
e(\mathbf{x}_i) := \Vert \mathbf{s}(\mathbf{x}_i)-\mathbf{f}_i\Vert_2, \qquad i=1,\ldots,n.
\end{equation*}
\end{enumerate}
\item While $\vert \{i:\, e(\mathbf{x}_i)>\epsilon\}\vert/n>\eta$ and $M \le M_{\mathrm{max}}$, repeat the following steps: 
\begin{enumerate}
\item ({\it marking}) for each $\ell=0,...,M-1$, mark the cells of level $\ell$ in ${\cal G_H}$ which are included in the support of a $B_J^{\ell}\,, J \in A^{\ell}$ such that: 
\begin{itemize}
\item there exists at least one data site ${\bf x}_i \in X \cap \mathrm{supp}(B_J^{\ell})$ such that $e(\mathbf{x}_i)>\epsilon;$  
\item 
there are at least $\lceil n_{\mathrm{loc}}/(n_1 n_2) \rceil$ data sites  in any of the $n_1 n_2$   subrectangles  which uniformly split $\supp(B_J^{\ell})$ (note that such splitting is just temporarily considered  to check whether this refinement requirement is satisfied);
\end{itemize}
\item ({\it update the hierarchical mesh}) update   ${\cal G_H}$ by dyadically  refining (in the two parametric directions) all the marked cells;
\item ({\it update the number of levels}) set $M$ equal to the current number of levels; 
\item  ({\it update the hierarchical space}) update  the THB-spline basis of $S_{\cal H}$ and the sets $A^{\ell}$, $\ell=0,...,M-1$;
\item ({\it spline update}) use Algorithm 1 to compute the coefficients of the hierarchical QI spline $\mathbf{s}$ defined in \eqref{hqi1} -- only coefficients associated   with THB-splines having new mother B-splines have to be computed; 
\item ({\it error update}) evaluate the new errors $e(\mathbf{x}_i), i=1,\ldots,n$ at the data sites;
\end{enumerate}
\end{enumerate}
\medskip
\noindent {\bf Outputs}
\begin{itemize} 
\item ${\cal G}_{\cal H}$: hierarchical mesh; 
\item  ${\cal T}({\cal G}_{{\cal H}}) $: THB-spline basis of $S_{\cal H}$;
\item coefficients of the hierarchical spline $\mathbf{s}$ defined in \eqref{hqi1}.
\end{itemize}

The refinement criterion has been motivated by the observation that, when the parameter values corresponding to the local data set $F_I^{\ell}$ are concentrated in a small part of the support of $B_I^{\ell}$, the quality of the approximation may be affected. For this reason, we consider a splitting of the two sides of $\supp(B_J^{M-1})$ in $n_1$ and $n_2$ uniform segments, respectively, and subdivide the support of $B_J^{M-1}$ in the resulting $n_1 n_2$ subregions, where we then check the presence of at least  $\lceil n_{\mathrm{loc}}/(n_1 n_2)\rceil$ data points. To simplify the usage of the algorithm by default  we set $n_1 = n_2 =1$  but different values can be chosen if suitable, see e.g. the data set with voids considered in Example~\ref{exm:3} of Section~\ref{sec:4}. 
 

Concerning the parameter $n_{\mathrm{loc}},$ we may observe that the requirement $n_{\mathrm{loc}} \ge n_{\mathrm{min}}$  guarantees that the points needed to compute the coefficients associated with the new functions in the first stage of  the next iteration can be found not too far from the support of the functions themselves. Indeed in the algorithm introduced in the previous section, after a few enlargements, $\Omega_J$ will surely include the support of a refined function of the previous level intersecting $\mathrm{supp}(B_J)$. As a consequence, analogously to $n_{\mathrm{min}}$, a high value of $n_{\mathrm{loc}}$ contributes to the reduction of oscillations deriving from overfitting, but this value should also be low enough to guarantee that the refinement strategy can generate a hierarchical spline space with enough degrees of freedom for satisfying the given tolerance $\epsilon$.  For this reason, some tuning is necessary for a good selection of $n_{\mathrm{loc}}$.

The proposed adaptive approximation method also extends to the case, not addressed in the previous works, of surfaces closed in one (or even two) parametric directions. Note that the local nature of the considered approximation approach makes the implementation especially easy, since coefficients associated with a THB present at successive steps of the adaptive refinement procedure (even if possibly further truncated) do not depend on such steps.


\section{Examples}\label{sec:4}
We present a selection of tests for the approximation of industrial data obtained by an optical scanning process of four different aircraft engine parts.  For each of these surfaces, as a characterizing dimension, we report the length $R$ of the diagonal of the minimal axis-aligned bounding box associated to the given point cloud. The parameter values are computed in all examples with standard parametrization methods based on a triangulation of the scattered data sets, see e. g., \cite{floater1997,floater2005}.  The bi-degree ${\bf d}$ is set equal to $(2,2)$ in the first considered example and always equal to $3$ in the other examples.
 The results highlight the effects of considering a minimum number of local data points (also) in the first stage of method, as well as the improvements obtained by introducing a regularized B-spline approximation for each local fitting with respect to the scattered data fitting scheme considered in \cite{bracco2018}. By combining these two changes, the two-stage approximation algorithm is more stable and unwanted oscillations are further reduced. 

Concerning the parameters in input to Algorithm 1, we have always set $\mu = 10^{-6},$ except for Example \ref{exm:2} where it was chosen even smaller.  In order to try to produce a very detailed reconstruction, a quite small value has been chosen for $n_{\mathrm{min}},$ always selecting it between $d^2$ and $(d+1)^2.$  Concerning the parameter selection for Algorithm 2, in all the presented experiments  we have set the maximum number of levels $M_{\mathrm{max}}$ equal to $8,$ fixing the percentage bound parameter $\eta$  equal to its default value ($\eta = 95 \%$). The error tolerance $\epsilon$ has been always set to $5 \cdot 10^{-5}$ m, except for Example \ref{exm:2} where we used the about halved value chosen for the same experiment in \cite{bracco2018}  (as a reference value, for each example consider the dimension $R$ characterizing the related point cloud). The integer parameters $n_1$ and $n_2$ also in input to Algorithm 2 have been always set to their unit default values, except for Example \ref{exm:3} which  required a different selection because of the voids present in the considered data set. The only parameter which has required a finer tuning for the reported experiments has been $n_{\mathrm{loc}}$ which is anyway always assumed greater than $n_{\mathrm{min}}.$  

\begin{example}\label{exm:1}
(\emph{Tensile}) In this example, we consider THB-spline approximations to reconstruct a part of a tensile from the set of 9281 scattered data shown in Figure~\ref{fig:exm1} (top) which has reference dimension $R = 2.5 \cdot 10^{-2}$ m. We compare the new local scheme based on local B-spline approximations with the algorithm based on local polynomial approximations of variable degree presented in \cite{bracco2018}, where this test was originally considered. Note that for this example we have never dealt with local sets of collinear points in our experiments.

As the first test, we ran both algorithms with the same setting considered in \cite{bracco2018}, namely, by starting with a $4\times 4$ tensor-product mesh with $\mathbf{d} =(2, 2)$, tolerance $\epsilon = 5\cdot 10^{-5}$ m,  percentage bound $\eta$ equal to the default $95 \%$ and $n_{\mathrm{loc}} = 20$. The algorithm with local polynomial approximations with the parameter choice considered in \cite{bracco2018} ($\sigma = 10^8$) led to an approximation with $2501$ degrees of freedom, $96.25\%$ of points below the tolerance and a maximum error of $1.22745 \cdot 10^{-4}$ m. The new scheme based on local B-spline approximations with $n_{\mathrm{min}}= 6$ and $\mu = 10^{-6}$  generated a THB-spline surface with $1855$ degrees of freedom that satisfies the required tolerance in $98.88\%$ of points with a maximum error of $8.06007\cdot 10^{-5}$ m. The number of levels used is $M = 5,$ but all the cells of the first two levels are refined.

As the second test, we ran both algorithms by starting with a $16\times 4$ tensor-product mesh with $\mathbf{d} =(2, 2)$, percentage $\eta$ equal to the default, tolerance $\epsilon = 5\cdot 10^{-5}$ m, and $n_{\mathrm{loc}} = 15$. The algorithm with local polynomial approximations led to an approximation with $5922$ degrees of freedom, $98.18\%$ of points below the tolerance and a maximum error of $1.44222 \cdot 10^{-4}$ m. The surface and the corresponding hierarchical mesh are shown in Figure~\ref{fig:exm1} (center). This approximation clearly shows strong oscillations on the boundary of the reconstructed surface, due to a lack of available data points for the local fitting in correspondence of high refinement levels. The scheme based on local B-spline approximations, with $n_{\mathrm{min}}= 7$ and $ \mu = 10^{-6}$   produced a THB-spline surface with $1960$ degrees of freedom that satisfies the required tolerance in the $99.36\%$ of the data points  with a maximum error of $8.10814\cdot 10^{-5}$ m. The surface, free of unwanted oscillations along the boundary, and the corresponding hierarchical mesh are shown in Figure~\ref{fig:exm1} (bottom). The number of levels used is $M = 4,$ with all the cells of level $0$ refined. 

\begin{figure}
\centerline{
\includegraphics[scale=.40]{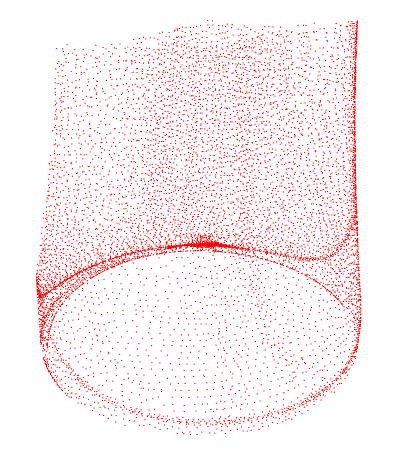}
\includegraphics[scale=.40]{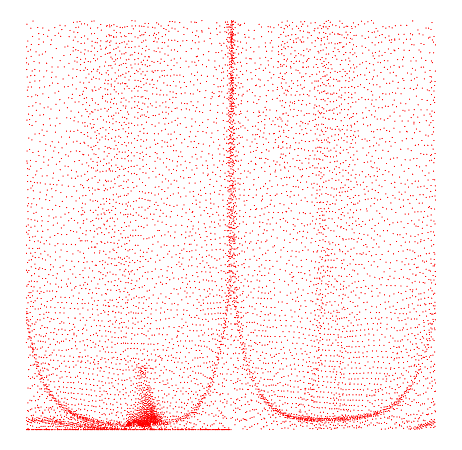}}
\centerline{\includegraphics[scale=.40]{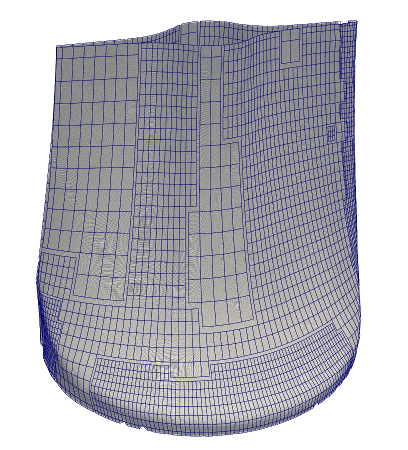}
\includegraphics[scale=.40]{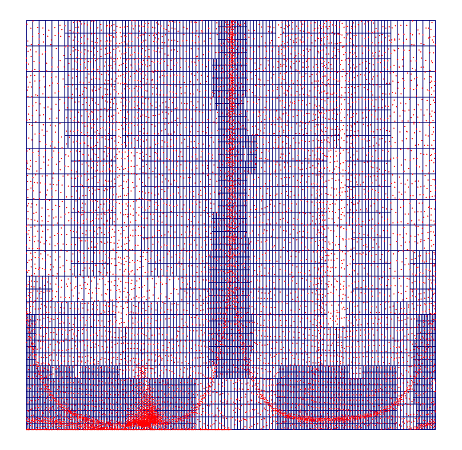}}
\centerline{\includegraphics[scale=.40]{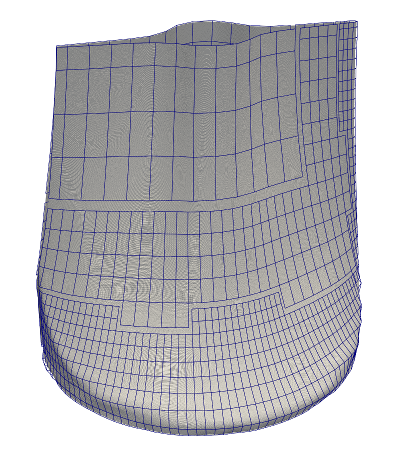}
\includegraphics[scale=.40]{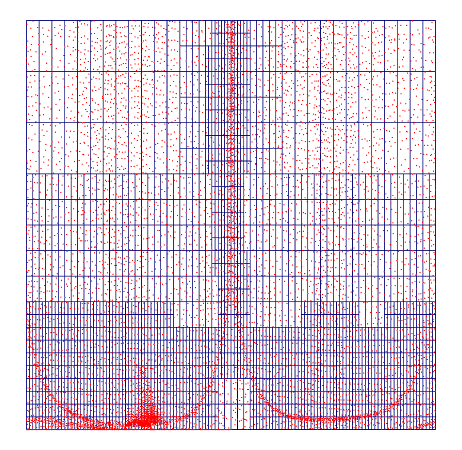}}
\caption{Example~\ref{exm:1}: scattered data set corresponding to a critical part of a tensile (top), the reconstructed surfaces and the corresponding hierarchical meshes obtained with the algorithm presented in \cite{bracco2018} (center) and the new scheme (bottom).}
\label{fig:exm1}     
\end{figure}

\end{example}

\begin{example}\label{exm:2}
(\emph{Blade}) In this example, we test the second example considered in \cite{bracco2018} on the set of 27191 scattered data representing a scanned part of a blade shown in Figure~\ref{fig:exm2} (top) whose reference dimension $R$ is equal to $5 \cdot 10^{-2}$ m. Again, to compare the new local scheme with the algorithm based on local polynomial approximations there considered, we ran both algorithms with the same setting of \cite{bracco2018}, namely, by starting with a $4\times 4$ tensor-product mesh with $\mathbf{d} =(3, 3)$, tolerance $\epsilon = 2\cdot 10^{-5}$ m,  percentage bound $\eta$ equal to the default $95 \%$  and $n_{\mathrm{loc}} = 60$. The algorithm with local polynomial approximations with the parameter choice considered in \cite{bracco2018} ($\sigma = 10^8$) led to an approximation with $12721$ degrees of freedom, $97.06\%$ of points below the tolerance and a maximum error of $1.08043 \cdot 10^{-4}$ m. The new scheme based on local B-spline approximations with $n_{\mathrm{min}}= 6$ and $\mu = 10^{-8}$  generated a THB-spline surface with $8314$ degrees of freedom that satisfies the required tolerance in $99.94\%$ of points with a maximum error of $1.32976\cdot 10^{-4}$ m. The surface and the corresponding hierarchical mesh are shown in Figure~\ref{fig:exm2} (bottom). The number of levels used with the new scheme is $7$ (the first two are not visible in the mesh because their cells are fully refined), one less than with the old approach. Finally for completeness we precise that the local collinearity check in Algorithm 1 is active only for 5 coefficients.
 
\begin{figure}[t!]
\centerline{\includegraphics[scale=.33]{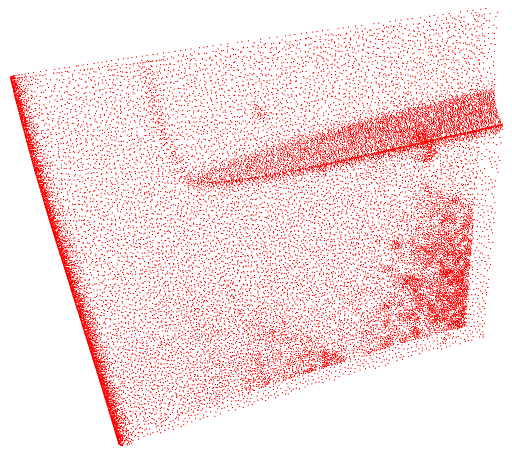}
\includegraphics[scale=.33]{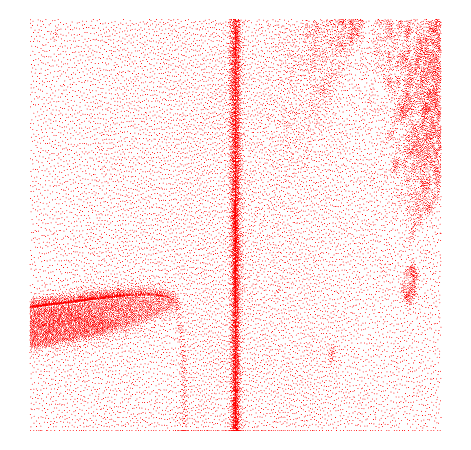}}
\centerline{\includegraphics[scale=.33]{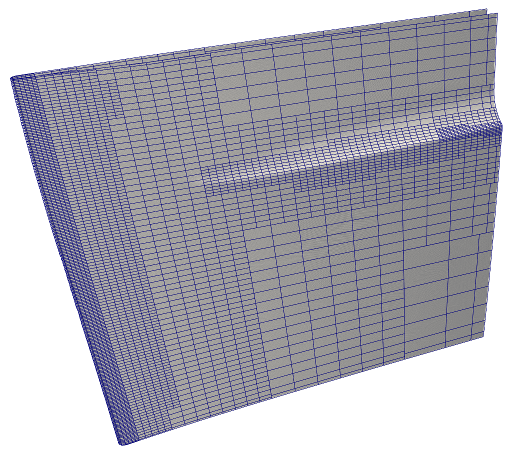}
\includegraphics[scale=.168]{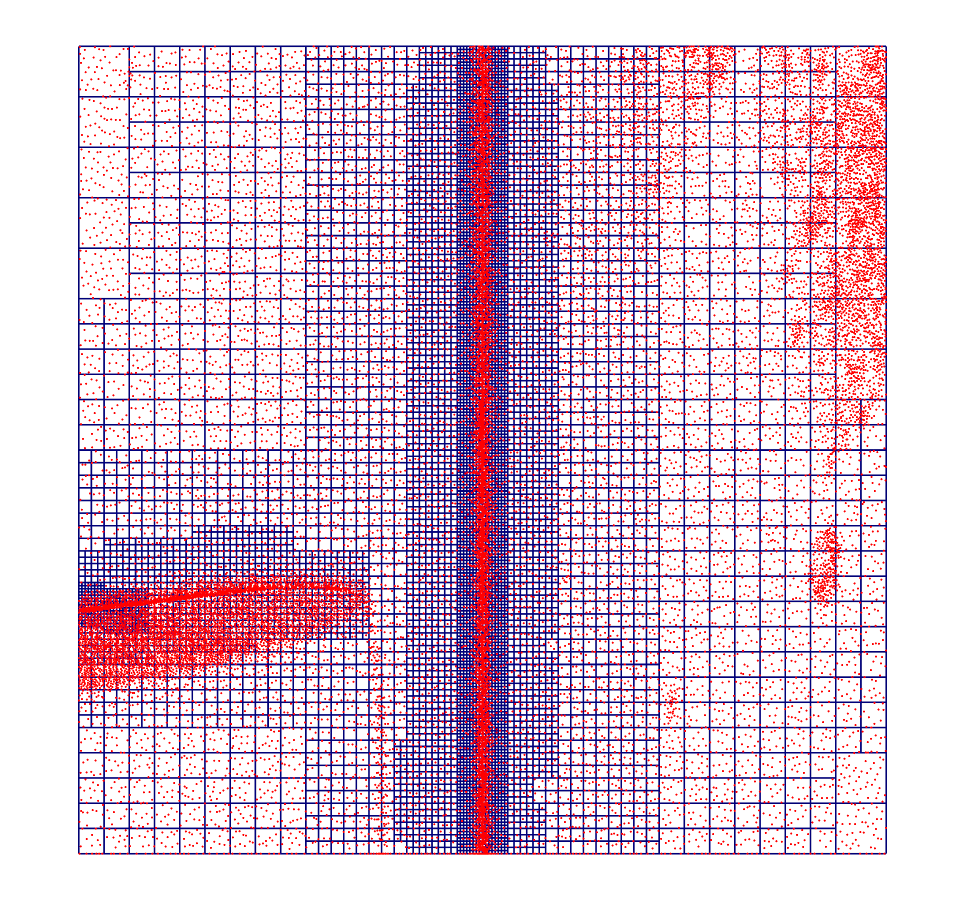}}
\caption{Example~\ref{exm:2}: scattered data set corresponding to a critical part of a blade (top), the reconstructed surface and the corresponding hierarchical mesh obtained with the new scheme (bottom).}
\label{fig:exm2}     
\end{figure}
\end{example}

\begin{example}\label{exm:3} (\emph{Endwall})
In this example we illustrate the behavior of the adaptive algorithm on data sets with voids by considering the reconstruction of an endwall part from the scattered data set of 43869 points shown in Figure~\ref{fig:exm3} (top) ($R = 5 \cdot 10^{-1}$ m). The figure shows that in this case the data set represents a model with three different holes, where no input data are available. The aim of this reconstruction is to avoid artifacts due to lack of points and obtain a sufficiently regular surface (e.g. by avoiding self-intersections), that can be post-processed with standard geometric software tools to obtain a suitably trimmed model. Consequently, not only the number of points in the local data sets is important to reach this aim, but also their distribution. To properly address this issue, we consider a real density parameter with value between 0 and 1 which determines whether the distribution of the points in the local set is reliable or not for the fitting. The distribution of the local data points is computed as the number of mesh cells of level $\ell$ inside the support of $B^\ell$ or its enlargement, which contain at least one point, over the total number of mesh cells, either in the support of $B^\ell$ or its enlargement. If this ratio is below a density threshold, then more data points are required and the function support is enlarged  for the computation of the local approximation in the first stage of the method. The approximation is developed by starting from a $32 \times 32$ tensor-product mesh, with $\mathbf{d}=(3,3)$, $n_{\mathrm{loc}} = 15,\, n_{\mathrm{min}} = 11$, $\mu = 10^{-6}$  and $n_1=n_2=2$.\ A choice of the density parameter $\delta$ equal to $0.3$ permits to take care of the difficult distribution of data points in the construction of the approximation. By considering a tolerance $\epsilon = 5\cdot 10^{-5}$ m and a percentage bound $\eta$ equal to the default $95\%$, the refinement  generated a THB-spline approximation with $M=3$ and $11211$ degrees of freedom, $98.70\%$ of points below the threshold and a maximum error of $5.68999 \cdot 10^{-4}$ m. The surface and the corresponding hierarchical mesh are shown in Figure~\ref{fig:exm3} (bottom). In this case there are $18$ coefficients of the last level and $15$ of the last but one (all associated with B-splines whose support intersects a void) such that the related $X_J$ is made up of all collinear points.

\begin{figure}[t!]
\centerline{\includegraphics[scale=.31]{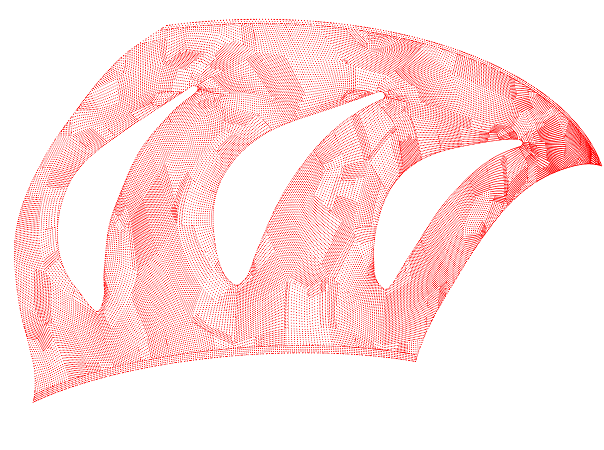}
\includegraphics[scale=.31]{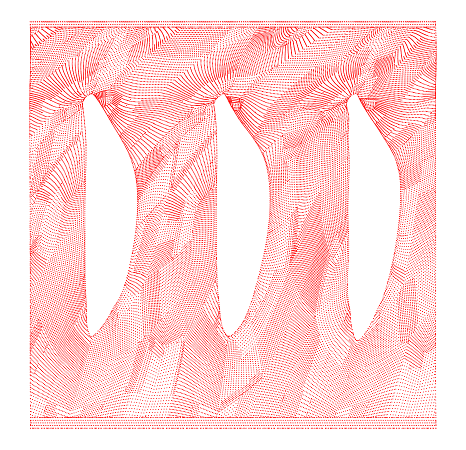}}
\centerline{\includegraphics[scale=.32]{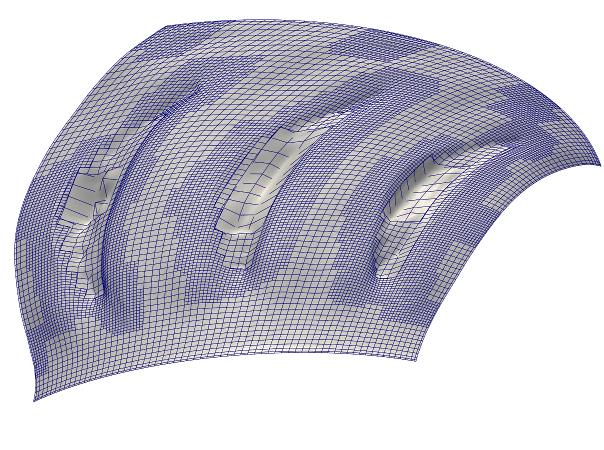}
\includegraphics[scale=.32]{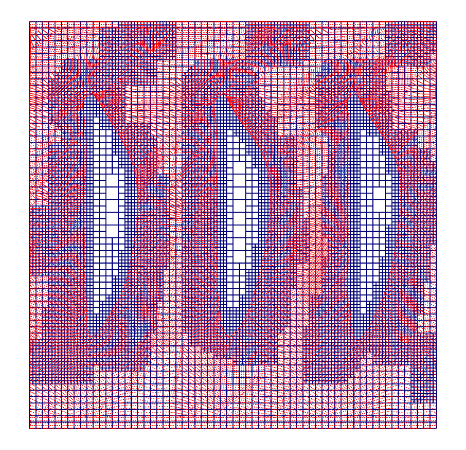}}
\caption{Example~\ref{exm:3}: scattered data set corresponding to a critical part of an endwall (top), the reconstructed surface and the corresponding hierarchical mesh obtained with the new scheme (bottom).}
\label{fig:exm3}     
\end{figure}

\end{example}

\begin{example}\label{exm:4}
(\emph{Airfoil})
This example illustrates the behavior of the new adaptive fitting algorithm with local B-spline approximations for  surfaces closed in one parametric direction. We test the scheme to reconstruct a blade airfoil from the set of 19669 scattered data shown in Figure~\ref{fig:exm4} (top) ($R = 10^{-1}$ m). We ran the method by starting with a $32\times 4$ tensor-product mesh with $\mathbf{d} =(3, 3)$,  setting  $\eta = 95 \%\,, \, \epsilon = 5\cdot 10^{-5}$ m, and $n_{\mathrm{loc}} = 30$ in Algorithm 2 and using Algorithm 1  with $n_{\mathrm{min}}= 12$ and $\mu = 10^{-6}.$   The refinement strategy produced an approximation with $M=3$ and $1856$ degrees of freedom distributed in the last two levels, that satisfies the required tolerance in $95.06\%$ of the data points with maximum error $1.87742 \cdot 10^{-4}$ m (observe also that, as well as for Example \ref{exm:1}, at the local stage $X_J$ is never made up of all collinear points). The surface and the corresponding hierarchical mesh are shown in Figure~\ref{fig:exm4} (bottom). By trying to force additional refinement, some oscillations appear. In this case, they are consistent with the data distribution since there are clusters of high density points, due to scan noise. Consequently, they do not represent artifacts caused by regions with very low density of data and cannot be prevented by exploiting the bound for cardinality of the local data sample governed by $n_{\mathrm{loc}}$.

\begin{figure}[t!]
\centerline{
\includegraphics[scale=.34]{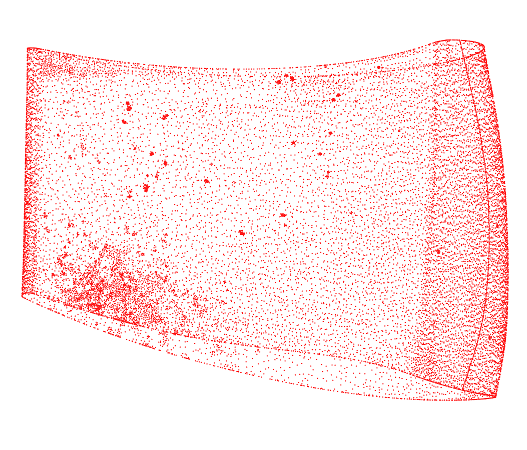}
\includegraphics[scale=.34]{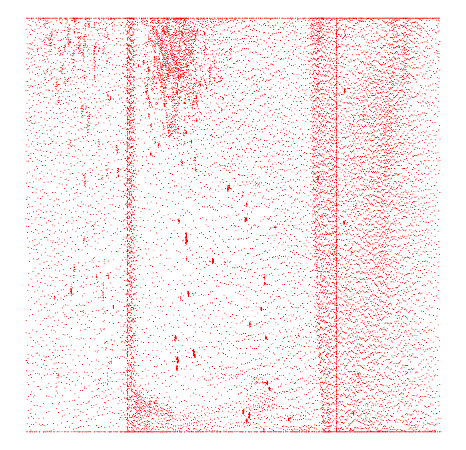}}
\centerline{\includegraphics[scale=.34]{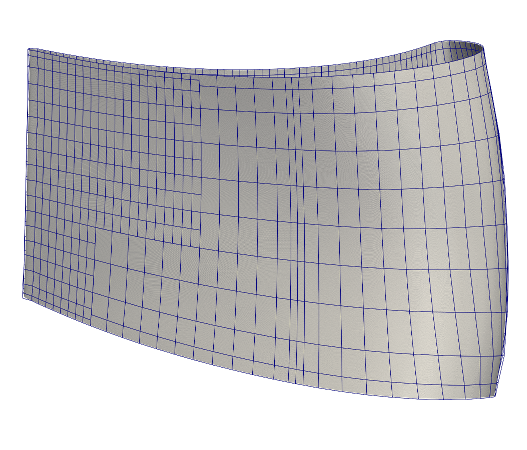}
\includegraphics[scale=.34]{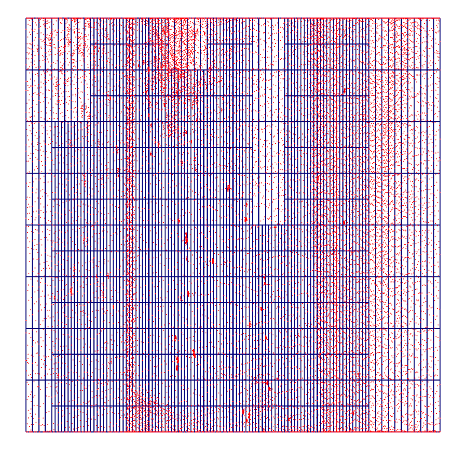}}
\caption{Example~\ref{exm:4}: scattered data set corresponding to a critical part of an airfoil (top), the reconstructed surface and the corresponding hierarchical mesh obtained with the new scheme (bottom).}
\label{fig:exm4}     
\end{figure}
\end{example}

 
\section{Appendix}
In this appendix we give the proof that, assuming ${\bf d} = (d_1,d_2)$ with $d_i \ge 1, i=1,2,$ the local vector spline $\mathbf{s}_J$ defined in Section \ref{sec:2} exists and is unique, provided that
the sites in $X_J$ are not collinear.
First of all we observe that  the objective function in (\ref{obj}) can be split in the sum of three analogous objective functions, one for each component $ s_J^{(k)}, k=1,2,3,$ of $\mathbf{s}_J,$   
$$
\sum_{i \in I_J} 
\Vert \mathbf{s}_J(\mathbf{x}_i) - \mathbf{f}_i \Vert_2^2 + \mu \, E(\mathbf{s}_J) = \sum_{k=1}^3  \bigl( \sum_{i \in I_j} \bigl( s_J^{(k)}(\mathbf{x}_i) - (\mathbf{f}_i)_k \bigr)^2   + \mu \, \rho( s _J^{(k)} ) \bigr) \,,$$
where
$$\rho \bigl( s _J^{(k)}\bigr) := \int_{\Omega_J} 
\left(\frac{\partial^2  s_J^{(k)}}{\partial x^2}\right)^2
+ 2 \left(\frac{\partial^2  s_J^{(k)}}{\partial x \partial y}\right) ^2
+\left(\frac{\partial^2  s_J^{(k)}}{\partial y^2}\right)^2 \,\text{d}x\text{d}y\,.$$
Thus the study can be developed in the scalar case and for brevity we remove the subscript or superscript $k$ ranging between $1$ and $3.$ 
The analysis is developed in the following theorem, where $s_J: \Omega_J \rightarrow \R$ denotes the local spline in ${\cal S}_J$ associated to the
coefficient vector $\mathbf{c} \in \R^{\ell_J},$ with $\ell_J:= | \Lambda_J |\,,$
$$s_J(\mathbf{x}) =  \sum_{I \in \Lambda_J}  c_I  B_I(\mathbf{x})\,.$$
\begin{thm}
{\it Let the considered spline bi-degree ${\bf d} = (d_1,d_2)$ be such that $d_i \ge 1, i=1,2.$
When the points $\mathbf{x}_i \in \Omega_J, i \in I_J$ are not collinear, there exists a unique local spline $s_J \in {\cal S}_J$ minimizing the following objective function,
\begin{equation} \label{obj2}
\sum_{i \in I_J} 
\bigl(s_J(\mathbf{x}_i) - f_i )^2 + \mu \, \rho(s_J)\,,
\end{equation}
where $\mu >0.$
If such points in $\Omega_J$ are collinear, then such minimizer does not exist or is not unique.}

\begin{proof} 
Let us observe that $\rho(s_J) = \mathbf{c}^T M \mathbf{c}\,,$ where $M \in \R^{\ell_J \times \ell_J}$ is such that 
$$M_{i,r} := \int_{\Omega_J} 
\left(\frac{\partial^2  B_I }{\partial x^2}\right)\left(\frac{\partial^2  B_R }{\partial x^2}\right)
+ 2 \left(\frac{\partial^2  B_I}{\partial x \partial y}\right) \left(\frac{\partial^2  B_R}{\partial x \partial y}\right)
+\left(\frac{\partial^2  B_I }{\partial y^2}\right)\left(\frac{\partial^2  B_R }{\partial y^2}\right)\,\text{d}x\text{d}y\,,$$
where we are assuming that, in the adopted ordering of the B-spline basis elements of ${\cal S}_J,$ $B_R$ and $B_I$ are the $r$--th and the $i$--th ones.
 On the other hand it is
 $$\sum_{i \in I_J} 
\bigl(s_J(\mathbf{x}_i) - f_i )^2  = \Vert V \mathbf{c} - {\bf F} \Vert_2^2\, = \mathbf{c}^T A^TA \mathbf{c} - 2{\bf F}^TA \mathbf{c} + {\bf F}^T{\bf F}\,,$$
where  ${\bf F} \in \R^{|I_J|}$ denotes the vector collecting all the $f_i, i \in I_J$ and $A$ is the $|I_J| \times \ell_J$ collocation matrix of the tensor-product B-spline basis generating ${\cal S}_J.$ 
Thus the objective function in (\ref{obj2}) can be written also as the following quadratic function,
$$\mathbf{c}^T (A^TA + \mu M) \mathbf{c}  - 2{\bf F}^TA \mathbf{c} + {\bf F}^T{\bf F}\,.$$
As well known a quadratic function admits a global unique minimizer if and only if the symmetric matrix  defining its homogeneous quadratic terms is positive definite and in such case the minimizer is given by its unique stationary point. In our case such matrix is $ A^TA + \mu M$ and the stationary points are  the solutions of the following linear system of $\ell_J$ equations in as many unknowns,
$$(A^TA + \mu M) \mathbf{c} = A^T {\bf F}\,.$$
Now, for all positive $\mu$ the matrix $A^TA + \mu M$ is symmetric and positive semidefinite since, for any vector $\bfz \in \R^{\ell_J}, \bfz \ne {\bf 0}$ it is $\bfz ^TA^TA \bfz \ge 0$ and $\bfz^TM\bfz \ge 0,$ the last inequality descending from the fact that $\bfz^TM\bfz = \rho(s_{\bfz}),$ where $s_{\bfz}(\mathbf{x}) = \sum_{I \in \Lambda_J}  \zeta_I  B_I(\mathbf{x})\,.$ Now if the points $\mathbf{x}_i, i \in I_j$ are distributed in $\Omega_J$ along a straight line $ax+by+c=0$, the proposition proved in Section~\ref{sec:2} implies that it is possible to find  $\bfz \in \R^{\ell_J}, \bfz \ne {\bf 0}$ such that $s_{\bfz}(\mathbf{x}) \equiv ax+by+c, \forall \mathbf{x} \in \Omega_J.$ This implies that $s_{\bfz}(\mathbf{x}_i )= 0\,,  \forall i \in I_J,$ that is the vector $A \bfz \in \R^{|I_J|}$ vanishes. On the other hand, clearly it is also $0 = \rho(s_{\bfz}) = \bfz^TM\bfz\,,$ since $s_{\bfz} |\Omega_J$ is a linear polynomial. This proves that the symmetric positive semidefinite  matrix $(A^TA + \mu M)$ is not positive definite when all the $\mathbf{x}_i, i \in I_J$ are collinear. This is the only possible data distribution associated to a non positive definite matrix. Indeed if the points $\mathbf{x}_i, i \in I_J$ are not collinear, if $\bfz \in \R^{\ell_J}, \bfz \ne {\bf 0}$ is associated to a non vanishing linear polynomial, it is $\bfz^TM\bfz = \rho(s_{\bfz})=0$ but $A \bfz \ne {\bf 0}$ and so  $\bfz^TA^TA\bfz > 0$; on the other hand if $\bfz \in \R^{\ell_J}, \bfz \ne {\bf 0}$ is not associated to a linear polynomial, then $\bfz ^TM \bfz = \rho(s_{\bfz})>0.$ 
\end{proof} 
\end{thm}

\section*{Acknowledgements}
We thank the anonymous reviewer for his/her useful suggestions. Cesare Bracco, Carlotta Giannelli and Alessandra Sestini are members of the INdAM Research group GNCS. The INdAM support through GNCS and Finanziamenti Premiali SUNRISE is also gratefully acknowledged.

\end{document}